\documentclass[12pt]{article}

\usepackage[utf8]{inputenc}       
\usepackage[T1]{fontenc}          
\usepackage{lmodern}              
\usepackage{geometry}             
\usepackage{graphicx}             
\usepackage{amsmath, amssymb}     
\usepackage[numbers]{natbib}      
\usepackage{hyperref}             
\usepackage{caption}              
\usepackage{pgfplots}
\usepackage{tikz}
\usepackage{algorithm}
\usepackage{algpseudocode}
\usepackage{fancyvrb}

\newcommand{\fitexpr}{\text{runtime} = a \times \text{NZ}^b}

\title{Empirical Asymptotic Runtime Analysis of Linear Programming Algorithms}
\author{Edward Rothberg \\
Gurobi Optimization, LLC \\
\texttt{rothberg@gurobi.com}}
\date{April 17, 2026}

\begin{document}

\maketitle
\begin{abstract}
  This paper takes an empirical look at asymptotic runtime growth
  rates for the most widely used algorithms for solving linear
  programming (LP) problems across a set of six optimization
  application areas that are known to produce large and difficult LP
  models.  On the algorithm side, we consider the simplex method,
  interior-point methods, and PDHG.  On the model side, we use a large
  language model (LLM) to create families of instances in different
  application areas, allowing us to study model types and sizes that
  are simultaneously synthetic and realistic.  The results indicate
  that simple regression models typically predict observed
  runtimes quite well within a model class, and
  that asymptotic behavior can vary significantly between the
  different algorithms.  This may have a significant impact on which
  algorithms will be most effective for solving large LP models in the
  future.
\end{abstract}

\section{Introduction}

Large linear programming (LP) models play a vital role in many
practical applications, including airline scheduling, supply-chain
planning, telecommunication network design, and electrical power
distribution.  The substantial time required to solve these problems
often limits the ability of these applications to solve real-world
problems in the desired timeframes.  As a result, substantial and
sustained research efforts, both academic and commercial, have been
devoted to advancing the state-of-the-art.  This led to the simplex
method~\cite{simplex51} in the 1950s, interior-point
methods~\cite{wright97} in the 1980s, and variants of the Primal-Dual
Hybrid Gradient (PDHG)
algorithm~\cite{pdlp22,pdhg11,hprlp24,cupdlp23,cupdlp+25} in the
2020s.  Each has been found to have strengths and weaknesses, with
none clearly dominating the others.

These algorithms have theoretical properties that can provide insights
into the tradeoffs between them.  Due to the variety and complexity of
practical optimization models, though, it is rare for pencil-and-paper
analysis to provide precise predictions of which algorithm would be
the best choice for a particular type of model.  Empirical
benchmarking is also important.

Extensive comparative benchmarking has in fact been performed between
implementations of LP algorithms.  To give one example, Hans
Mittelmann maintains a benchmark page~\cite{hans25} that compares
specific implementations of different algorithms for sets of LP
instances.  Such efforts are important for understanding which methods
are likely to perform best on which models, and also to track progress
over time.

Such benchmarking efforts have a significant limitation, though: they
are point-in-time studies.  They capture results for a particular set
of models on a particular set of machines.  One constant in linear
programming is that improved solver capabilities lead to larger models
(which lead to new improvements, which lead to...).  It is tough to
use point-in-time results to predict what will happen over time as
model sizes continue to grow.

Looking at runtime growth rates requires sets of optimization model
instances that are similar in character but different in size.  In
this paper, we use an LLM (ChatGPT 5.4 Thinking~\cite{chatgpt2026}) to
produce synthetic model generators.  Specifically, the LLM generates
optimization models from six different important application domains,
producing a family of models of varying sizes for each, where the
smallest is straightforward to solve and the largest is challenging.
We solve these model families using different algorithms, and use a
regression model to predict runtime from model size for each model
family and algorithm.  The goal is to gain insight into likely future
runtimes as model sizes continue to grow.

\section{The Algorithms}

A linear programming problem can be stated in its primal and dual
forms as:
\begin{equation*}
\begin{minipage}{0.45\linewidth}
\begin{align*}
\min \quad & c^\top x \\
\text{s.t.} \quad & A x = b \\
& x \ge 0
\end{align*}
\end{minipage}
\hspace{0.2cm}
\begin{minipage}{0.45\linewidth}
\begin{align*}
\max \quad & b^\top y \\
\text{s.t.} \quad & A^\top y + z = c \\
                  & z \ge 0
\end{align*}
\end{minipage}
\end{equation*}
where $A$ is an $m \times n$ {\em constraint matrix\/}, $c$ is the
{\em objective vector\/}, $b$ is the {\em right-hand-side vector\/},
$x$ is the {\em primal solution\/}, $y$ is the {\em dual solution\/},
and $z$ is the {\em reduced-cost vector\/}.

Any trio of vectors $(x,y,z)$ has associated primal and dual residual
vectors:
\begin{align*}
r_P &= b - Ax \\
r_D &= c - A^\top y - z
\end{align*}
A solution is optimal if (i) the primal and dual residual vectors are
zero, (ii) $x$ and $z$ are non-negative, and (iii) the objective gap is
zero. This last condition can be expressed directly as
$c^\top x = b^\top y$, or somewhat less directly as $x^\top z = 0$
(referred to as {\em complementarity\/}).

The three main algorithms available for solving LP problems are the
simplex method, interior-point methods, and PDHG.  In typical
situations, these are supplemented with a few additional algorithms,
including presolve, fill-reducing ordering, and crossover.  The mix
depends on the method.

\subsection{Presolve}

The first step in solving nearly any LP is
presolve~\cite{presolve1995}, which attempts to reduce the size of the
original LP model in order to reduce the cost of the subsequent solve.
Presolve is typically quick, but it has substantial sequential
components, so the cost can sometimes be significant relative to the
cost of a highly parallel LP algorithm.  Experience so far has
indicated that, even for the most parallel of LP algorithms, an
initial presolve step is still quite beneficial.  It reduces overall
solution time, and it frees modelers from the chore of having to
remove obvious redundancies from their models.

Once presolve is complete, the presolved model is passed to an LP
algorithm.  Once a solution is found, it is mapped back to a solution
of the original model.

\subsection{Simplex Method}

The simplex method~\cite{simplex51} is the oldest and most versatile
of the available algorithms for solving LPs.  It maintains a basis at
all times, which is a set $B$ of columns of size $m$ (the number of
constraints in $A$).  The associated solution is computed from the
submatrix induced by $B$: $x_B = A_B^{-1} b$, $x_N = 0$ (where $N$ is
the set of non-basic columns), and $y = A_B^{-T} c_B$.  In each
iteration (or {\em pivot\/}), a basic variable and a non-basic
variable swap places, giving a new $x_B$.  Swaps are chosen to improve
the objective function value, and iterations proceed until no
improving candidates are available, at which point the current iterate
is declared optimal.

Of course, small violations of LP optimality conditions are allowed in
the returned solution.  Typical termination criteria in simplex are
$\|r_P\|_{\infty} \leq \varepsilon$ and
$\|r_D\|_{\infty} \leq \varepsilon$, with $\varepsilon = 10^{-6}$
being the traditional choice.

The simplex method has primal and dual variants, where the dual is
effectively just the primal algorithm run on the dual problem.
Despite this symmetry, the dual version is much more effective in
practice.

A basis is a powerful construct in LP.  Its many virtues include: (i)
it produces highly accurate solutions, (ii) it represents solutions
quite concisely, and (iii) it enables quick reoptimization after
small changes to the problem (when used in conjunction with the
simplex method).

One significant limitation of the simplex method is its lack of
parallelism, which makes it ill-suited to today's multi-core
processors.  Simplex requires the solution of an extremely sparse
linear system in each iteration, and that system change slightly in
each iteration, making parallelization quite difficult.  A typical LP
solve requires thousands or millions of pivots, and one pivot depends
on the results of the previous one, making it difficult to parallelize
across iterations.

\subsection{Interior-Point Method}

Interior-point solvers~\cite{wright97} also maintain a current
iterate, but the process of improving that iterate is quite different.
The dominant computation in an interior-point solve is finding an
improving direction, which requires the solution of a weighted normal
equation $ADA^\top d = b$, where $A$ is the constraint matrix, $D$ is
a diagonal matrix, and the resulting $d$ vector is used to derive an
improvement direction.  The solution of this linear system is
typically performed using a sparse Cholesky factorization of
$ADA^\top$, followed by a pair of linear solves with the resulting
triangular factor matrix.  Iterations proceed until the primal
residual, dual residual, and complementarity are small enough.

Interior-point solvers typically aim to adhere to the same feasibility
tolerances as simplex, but will settle for looser tolerances if the
simplex tolerances appear unattainable.  Specifically, an interior-point
solver will typically fall back to
$\|r_P\|_2 \leq \varepsilon (1 + \|b\|_2)$
and $\|r_D\|_2 \leq \varepsilon (1 + \|c\|_2)$,
with $\varepsilon = 10^{-8}$ typical.  Termination also requires small
complementarity, so requiring $x^\top z / n \leq \varepsilon$
is also typical.

Sparse Cholesky factorization is quite amenable to parallelization, so
interior-point methods can make excellent use of multi-core computer
architectures.  Interior-point methods require many other steps,
though, several of which are not nearly as well suited to
parallelization.

One such step in interior-point solvers that should be highlighted is
{\em fill-reducing ordering\/}~\cite{georgeliu81}, which is an
integral part of sparse Cholesky factorization.  The cost of
performing sparse factorization depends heavily on the ordering of the
rows and columns of the sparse matrix being factored, so a heuristic
reordering is essential for reducing this cost.  Ordering only needs
to be performed once, since the structure of the sparse matrix being
factored remains fixed.  This step can still be expensive,
especially when running on machines with lots of cores, since the
factorization parallelizes well and ordering algorithms in general do
not.

Interior-point solvers are {\em locally quadratically convergent\/},
so once the iterate is close to the optimal solution,
accuracy improves quite quickly with each additional iteration. While
interior-point solvers typically produce extremely accurate solutions,
the resulting solutions do not share many of the nice properties that
a basic solution provides.  Thus it is typical for an interior-point
solve to be followed by a {\em crossover\/}~\cite{crossover91} step,
where the interior solution is transformed into a basic solution.  We
will say more about crossover shortly.

\subsection{Primal-Dual Hybrid Gradient Method}

The Primal-Dual Hybrid Gradient (PDHG)
algorithm~\cite{pdlp22,pdhg11,hprlp24,cupdlp23,cupdlp+25} is a
new method for solving LPs.  It also maintains a current iterate and
computes direction vectors to continually improve that iterate.  One
very practical difference versus interior-point methods is that the
direction vector computation is quite simple; it requires a pair of
sparse matrix-vector multiplications, a few dense vector additions and
dot products, and some scalar operations.  These operations map very
well to GPU architectures, making the algorithm quite competitive with
simplex and interior-point when run on a modern GPU.  PDHG is a
first-order method, so while each iteration is inexpensive,
convergence can be quite slow.  As a result, the method is typically
terminated with a solution that is much less accurate than those
produced by the alternative methods.

PDHG solvers traditionally terminate when
$\|r_P\|_2 \leq \varepsilon (1 + \|b\|_2)$,
$\|r_D\|_2 \leq \varepsilon (1 + \|c\|_2)$, and
$|c^\top x - b^\top y| \leq \varepsilon (1 + |c^\top x| + |b^\top
y|)$.  A value of $\varepsilon = 10^{-4}$ is common in academic
papers, but commercial implementations appear to have settled on
$10^{-6}$ as the standard default value.  Note that these termination
criteria are much looser than those for the simplex and interior-point
methods.

PDHG can also be followed by a crossover step to produce a basic
solution.  One obvious appeal is that crossover turns a potentially
low-accuracy starting solution into a highly accurate basic solution.
One potential downside is that a less accurate starting solution can
seriously harm crossover performance.

\subsection{Crossover}

Any primal/dual optimal solution pair, no matter how it was computed,
can be transformed to a basic optimal solution using a crossover
procedure~\cite{crossover91}.  The crossover algorithm is very similar
to simplex, with nearly the same parallelism limitations.  The
algorithm can be shown to run in polynomial time and is generally
quick in practice.

Theoretical results about crossover behavior rely on an assumption
that the starting solution and all intermediate solutions within the
algorithm achieve exact primal feasibility, dual feasibility, and
complementarity.  This is almost never the case in practice, due to
noise in the initial solution or accumulation of numerical errors
during the method.

It is well known that the further you stray from this assumption, the
more likely it is that crossover will run into trouble.  Two practical
consequences of this observation for this paper are: (i) we only
launch crossover if the preceding interior-point or PDHG solve
achieves its convergence criteria, and (ii) we tighten the
interior-point convergence criteria to reduce errors in crossover's
starting solution.  While tightening tolerances may sound expensive,
the local quadratic convergence behavior of interior-point methods
means that this usually only adds a few iterations.  While we could also
tighten the convergence criteria for PDHG, as a first-order method
that will typically lead to significant increases in PDHG runtimes.

What does it look like for crossover to run into trouble?  The typical
failure mode is for crossover to finish with a basis with substantial
primal and/or dual infeasibilities.  This basis is then used to
warm-start the general-purpose simplex method, which can require a
substantial amount of time.

\subsection{Accuracy}

Solution accuracy has come up a number of times in this discussion.
It may be helpful to quantify the differences for the methods
considered here, using results from another recent
paper~\cite{pdhghybrid26}.  Specifically, that paper measured the
maximum absolute violation of any primal or dual constraint for a set
of models, and then compared geometric means of this value over a
benchmark set for several LP algorithms.  For the Mittelmann LP
set~\cite{hans25}, which consists of 43 models of modest size and
difficulty, using Gurobi default tolerances, simplex gave a mean
absolute constraint violation of roughly $10^{-10}$, interior-point
without crossover gave a mean violation of roughly $10^{-8}$, and PDHG
without crossover gave a mean violation of roughly $10^{-4}$.  One can
certainly debate how much accuracy is truly necessary, but the
differences between the methods can be quite stark.

\section{The Model Generator}

As noted earlier, we chose to study optimization models from six
application areas that are known to require the solution of large,
challenging linear programming problems.  We asked an LLM to create a
synthetic model generator for each application area and then used that
generator to create the instance families that formed our test sets.
This section provides details on that process.

\subsection{Model Types}

The first task we undertook was to identify a set of high-impact
application areas that depend heavily on solving difficult linear
programming problems.  In most cases, these LPs are actually
relaxations of mixed-integer programming (MIP) problems, but solving
relaxations is almost always a vital first step in solving the
associated MIP.

The application areas we chose were:
\begin{itemize}
\item Airline Fleet Assignment (\textbf{Fleet}): An airline fleet
  assignment model chooses specific aircraft to fly specific legs of a
  schedule in order to best match aircraft capacity to expected demand.
  Our synthetic model captures a number of characteristics of
  the flight network, including the number of aircraft, airports, and
  flights, the number of maintenance facilities, necessary maintenance
  intervals, and the time horizon. Integrality restrictions
  are ignored, giving a linear relaxation.
\item Gas Network Optimization (\textbf{GasNet}): Gas transmission
  network optimization optimizes the distribution of natural gas over
  a pipeline network, using pressure generated by compressors to move
  gas from producers to consumers over some time
  horizon while minimizing operating costs.  The model captures the
  number of supply and demand nodes, the number of transit nodes, the
  number of compressors, and the number of time periods.  The model
  performs a linear approximation of the non-linear physics that
  controls the behavior of a network under pressure.  It also
  relaxes integrality restrictions.
\item Production Planning (\textbf{ProdPlan}): A production planning
  model decides how to schedule the production and distribution of
  finished goods in order to effectively use available resources to
  satisfy customer orders.  The model captures the number of plants,
  distribution centers, products, and product families, as well as
  inventories and storage capacities for the various inputs and
  outputs in different locations, optimized over a time horizon.
  Integrality restrictions are ignored, giving a linear relaxation.
\item Supply-Chain Network Design (\textbf{SCND}): A supply-chain
  network design problem determines the best configuration of a supply
  network given typical customer demand patterns.  It chooses
  distribution center locations and sizes, and the associated shipping
  lanes to minimize costs, minimize late deliveries, and meet CO2
  emissions targets.  The model captures: (i) supplier data: locations,
  maximum supply rates, and quantities, (ii) customer data:
  locations and demands, (iii) distribution center data: capacities
  and throughputs, (iv) transportation lanes: inbound and outbound
  capacity.  Integrality restrictions are ignored, giving a linear
  relaxation.
\item Telecommunication Network Design (\textbf{TelecomND}): A
  capacity planning model that decides how to configure a
  telecommunication network backbone to satisfy expected demand over
  multiple commodities and time periods in the most efficient way.
  The model captures the topology of the network (number of clustered
  POP facilities, communication nodes and links, link capacities),
  demands placed on the network over time (by
  commodity), and service-level agreements (latency tolerances).
  Integrality restrictions are ignored, giving a linear relaxation.
\item Electrical Power Unit Commitment (\textbf{UnitCommit}): An
  electrical power generation model creates a market for electricity,
  matching power producers with expected power consumers over a time
  horizon at agreed-upon prices.  The model captures generator
  information (number of generators, startup and shutdown costs, wind
  and solar production), demand information (expected hourly demand
  and reserve requirements), and storage management information
  (available storage capacity, charge/discharge choices).  Integrality
  restrictions are ignored, giving a linear relaxation.
\end{itemize}
Later results will refer to these test sets using the abbreviated
names in parentheses above.

\subsection{Model Generator}

Due to recent rapid improvements in Generative AI technology, large
language models (LLMs) are now quite skilled at writing and debugging
code.  While their competency depends on the amount of relevant
training data available, recent experience suggests that LLMs are now
quite skilled at generating code for many specialized domains,
including optimization.  Domain-specific tools like
Gurobot~\cite{gurobot} and OptiMind~\cite{zhang2025optimind} aim to
supplement these general-purpose tools by adding specialized domain
knowledge, but tools like ChatGPT~\cite{chatgpt2026} are also quite
capable of building optimization models.

In this paper, model generators are created by providing the prompt
shown in Figure~\ref{fig:mainprompt} to {\em ChatGPT 5.4 Thinking\/}.
\begin{figure}[ht]
\centering
\footnotesize
\begin{minipage}{0.8\linewidth}
\begin{Verbatim}[frame=single]
[Description of model class goes here]

Create a linear programming model for this problem.  The maximum
problem size should be quite large - tens of millions of constraints
and variables.  The intent is to use this model in a benchmarking
context.

The model should be as realistic as possible.  One aspect of this is
that the input parameters should not substantially exceed realistic
real-world values.  For example, don't suggest an airline fleet
assignment model with more than 2000 aircraft; that's more than any
real airline has.  To give another example, don't consider a ten-year
planning horizon in a production planning model - realistic horizons
are at most a few years.

Be careful to only generate models that have feasible solutions.

The generator should call the gurobipy API, and the output should be a
.lp file.  The name of the .lp file should capture the input parameter
values used to generate the instance.

Take your time creating the model.  In particular, do research to
understand the sorts of variables, constraints, and objective
functions that are typical in optimization models of this type.

Propose a set of input parameter values that would make a nice
'benchmark ladder', producing a family of benchmark models of
different sizes and difficulties.
\end{Verbatim}
\end{minipage}
\caption{Main LLM prompt}
\label{fig:mainprompt}
\end{figure}
The missing piece in this prompt is the description of the
optimization problem to be solved.  We found that these descriptions
could be quite high-level, simply providing the name of the problem and a
list of features we wanted the model to capture (e.g., respecting
maintenance requirements for a fleet assignment model).  In fact, the
application descriptions from the previous section could be used for
this purpose.  For three of our six model families, these instructions
alone were sufficient to obtain working models.  In the other cases, a
bit of iteration was required: twice because the resulting models were
infeasible and once because the model generator was too inefficient to
generate large models in a reasonable amount of time.

Instance sizes were chosen by evaluating points on the suggested
``benchmark ladder'', where a specific instance is defined by a set of
values for the various input parameters to the model generator.  At
one extreme, we looked for instances that
could be solved within around 100 seconds.  At the other extreme, we
looked for instances that two of the LP algorithms could solve within
10,000 seconds.  We generated 100 instances for each model class by
linearly interpolating input parameters between the smallest and
largest instance sizes.  The intent was to obtain a set of models that
spans a wide range of sizes while at the same time allowing us to
compare the different algorithms without requiring an inordinate
amount of computing time to run the full benchmark set.
Table~\ref{tab:model_sizes} shows the sizes of the largest models
in each set.
\begin{table}[htbp]
\centering
\begin{tabular}{l|rrr}
           & Constraints (M) & Variables (M) & Non-zeros (M) \\ \hline
Fleet      & 1.9 & 32.1 & 71.0 \\
GasNet     & 5.2 & 5.5 & 21.9 \\
ProdPlan   & 8.8 & 35.9 & 82.4 \\
SCND       & 2.0 & 21.5 & 150.1 \\
TelecomND  & 2.1 & 17.7 & 109.6 \\
UnitCommit & 8.2 & 7.0 & 53.8 \\
\end{tabular}
\caption{Size of the largest instance of each model type (after presolve).}
\label{tab:model_sizes}
\end{table}

\subsection{Realism}

How realistic were the models that the LLM produced?  An optimization
model instance is built from two inputs: the abstract model and the
data.  The abstract model describes the overall structure of the
model: the variables whose values must be chosen, the objective
function on these variables, the types of constraints that need to be
captured, etc.  Given the high-level nature of a typical model, it is
quite feasible to evaluate realism.  We are familiar with the problems
being solved here, and a quick perusal of the code generated by the
LLM suggested that the model had correctly captured the features of
the problem that it claimed to be capturing.

Input data realism is much harder to gauge.  While an abstract model
can be fully understood from a brief textual description, the data for
a large model involves potentially millions of values, and the overall
behavior of the model could depend heavily on each one.

Generation and use of synthetic data is an active research topic in
the LLM community~\cite{synthetic25}, with some notable successes and
failures.  The difference between this work and ours is that outcomes
can be easily evaluated when using sythetic data to train an LLM:
train the LLM with and without the synthetic data and compare the
resulting models' predictive power.  In our context, it is difficult
to propose an objective measure for optimization model realism.

While we are confident that our LLM-generated optimization models are
reasonable facsimiles of realistic abstract optimization models, we
are much less confident that the synthetic data that was used to
generate specific instances reflects realistic business scenarios.  We
are hopeful that this general approach can provide insight into the
behavior of practical models, but ultimately they are still synthetic
models.

\subsection{Gurobi Model Library}

An obvious question is why rely on synthetic benchmarks at all when
Gurobi has a vast library of benchmark models, collected from
customers and academic users over nearly two decades.  There
are a few reasons for this.

The first is that a large fraction of the models in our library were
sent to us in a technical support context; a customer was having
trouble solving a model, so they sent it to us to get help in
understanding and fixing an issue.  That means that a significant
fraction of the models in our test set are pathological in some sense.
While such models are important for improving product robustness, they
are less well suited for evaluating performance for common model
instances.

Another reason is that really large models can be unwieldy: large
file sizes, large runtimes, etc.  As a result, we typically do not get
that many.  If someone wants help with a model, their aim is typically
to send us the smallest model that illustrates the issue they are
trying to address.

Finally, it is rare for us to receive a large family of related
models.  To predict runtimes from model sizes, we really want models
of many different sizes (we generated 100 for this paper).  While you
could imagine collecting multiple model instances from different
sources that are from the same application area, the reality is that
different people will generate different models for the same problem,
and these different models can have very different performance
characteristics.

\section{The Regression Model}

Our basic approach to exploring the asymptotic behavior of the various
LP algorithms on our benchmark models is to run a variety of problem
sizes and use regression to compute an expression that explains the
resulting runtimes in terms of the size of the instance.  We greatly
simplify the task by assuming that runtime will be driven by a
``power-law'' function ($\text{runtime} = a \times \text{size}^b$),
and that we can capture the size of the instance using just the number
of non-zeros in the constraint matrix.

On the second point, we have verified that the non-zero count in $A$
is almost exactly a polynomial function of the number of rows and
columns in $A$ for all of our model families.  Including dependent
variables in regression is known to cause spurious results, so
building an expression on a single independent variable is likely to
produce more consistent results.

Regarding the simple power-law regression fit, our hope is that we're
solving large enough instances of these models that runtime will be
dominated by the term with the largest exponent.  The advantage of a
power-law formulation is that the $a$ and $b$ terms can be estimated
using log-log linear regression - we can use ordinary least squares
(OLS) to compute the best fit between the log of the non-zero count
and the log of the runtime.

Acknowledging that our goal of fitting the data with a power-law may
overly optimistic, we also look for evidence that a better fit might
be possible.  To do so, we use non-linear regression~\cite{regression}
to fit polynomial, exponential, and logarithmic functions to the data.
When we notice low quality in the power-law regression fit, we compare
against the best fit found with any of the alternatives.

To estimate the predictive performance of the regression model, we
report 95\% confidence intervals and
the $R^2$ measure~\cite{regression} (also known as the {\em coefficient of
  determination\/}).  The latter is a standard measure of the extent to which
variations in the results can be explained by the predictive model.
When it comes to setting a threshold for $R^2$ in order to declare the
model to be ``good'', there are no universal answers; it depends on
what you plan to do with the results.  Based on guidelines from other
problem domains and observations of our results, we interpret $R^2$
values above $0.7$ (in log space) as being a sign of an effective model.  Boiling
data points down to a single number can have limitations, of course,
so we will also sometimes include plots of the data points to
understand cases where the predictive power is poor.  In most cases,
poor quality models are due to noise in the data rather than
limitations of the function we are trying to fit.

We should add that the predicted exponent in the power-law model
becomes less reliable as the $R^2$ value drops, so these numbers
should not be taken too seriously when $R^2$ is below $0.7$.

One objection that can reasonably be raised is that the runtimes of LP
algorithms depend on two factors, the number of iterations required
and the cost of each iteration.  Iteration cost can reasonably be
expected to be a function of the size of $A$ (a linear relationship
for PDHG, a more complex function involving graph separator sizes for
interior-point methods).  The number of iterations is more complex,
often depending on numerical properties of $A$, and thus perhaps not
as easily captured as a function of the size of $A$.  The number of
iterations required by PDHG, for example, is known to depend on the
condition number of $A$~\cite{pdhg11}.  We should note that
conditioning for a class of matrices often bears a direct relationship
to the size of the matrix.  To give a simple example, the condition
number of a 2-dimensional $N \times N$ Laplacian matrix is known to be
proportional to $N^2$~\cite{laplacian}.  While the constraint matrices
considered here are too complex to readily admit closed-form condition
number expressions, it seems quite plausible that conditioning will be
related to the size of the constraint matrix in some simple way.

\section{Experimental Results}

Before presenting measured results, let us first describe our testing
environment, including details on the software and machines used.

\subsection{Testing Environment}

All of our tests were performed using Gurobi version 13.0, which
includes state-of-the-art implementations of simplex, interior-point,
PDHG, and crossover algorithms.  Parameters were left at default
values with the following exceptions: (i) the {\em Method\/} parameter
was changed to switch between LP algorithms, (ii) the {\em PDHGGPU\/}
parameter was set to 1 when running on a GPU, (iii) the {\em
  BarConvTol\/} parameter was set to $10^{-12}$ to tighten the
interior-point convergence criteria and improve the robustness of the
crossover step.  Due to the interior-point algorithm's locally
quadratic convergence behavior, the cost of the last change is
typically only a small number additional iterations.

Most of our computational tests were performed on
a large cluster of AMD EPYC 4364P-based systems, each with an 8-core
processor and 128~GB of memory.  As CPU-only systems with limited
memory bandwidth, these machines give quite poor performance for
PDHG.  We use a 10,000-second time limit for our simplex
and interior-point tests, and a 50,000-second time limit
for our PDHG tests.

We also perform some tests on a GPU system that contains an AMD EPYC
9575F CPU and an Nvidia RTX Pro 6000 GPU.  The CPU has 64 cores and
768~GB of memory, and the GPU has 96~GB of HBM memory.  The CPU and
GPU in this system are both from the latest generation, although
neither is top-of-the-line.  For large models, the CPU is roughly 8
times faster than the EPYC 4364P for interior-point iterations, and
the GPU is roughly 40 times faster than the EPYC 4364P CPU for PDHG
iterations.

\subsection{Results}

Recall from our earlier discussion that the various algorithms
for solving LP models consist of a number of building blocks:
\begin{itemize}
\item \textbf{Simplex}: presolve, simplex iterations
\item \textbf{Interior-point}: presolve, fill-reducing ordering, interior-point iterations, crossover
\item \textbf{PDHG}: presolve, PDHG iterations, crossover
\end{itemize}
Our asymptotic analysis will look at these pieces
individually to try to gauge which is likely to dominate
runtime as model sizes grow.

Let us first look at predicted runtime growth rates, using a power-law regression fit,
for the algorithm iterations alone - excluding presolve, crossover, etc.
All of these results were computed using our CPU-only EPYC 4364P system.
These numbers are shown in Table~\ref{tab:iteration_fit}
(and also in Figure~\ref{fig:iteration_fit}).
\begin{table}[ht]
\centering
\begin{tabular}{l|cc|cc|cc}
  \textbf{} &
   \multicolumn{2}{c|}{Dual simplex} &
   \multicolumn{2}{c|}{Interior-point} &
   \multicolumn{2}{c}{PDHG} \\
    &
   Exponent (CI) & $R^2$ &
   Exponent (CI) & $R^2$ &
   Exponent (CI) & $R^2$  \\ \hline
Fleet      & -           & -    & [1.68,1.92] & 0.92 & - & - \\
GasNet     & -           & -    & [2.62,2.82] & 0.97 & [0.56,1.14] & 0.26 \\
ProdPlan   & [1.46,1.64] & 0.93 & [1.64,1.71] & 0.99 & [1.03,1.07] & 0.99 \\
SCND       & -           & -    & [1.58,1.66] & 0.99 & [1.24,1.36] & 0.95 \\
TelecomND  & [1.67,2.27] & 0.81 & [1.48,1.55] & 0.99 & [1.26,1.44] & 0.91 \\
UnitCommit & [1.92,2.06] & 0.97 & [1.57,1.71] & 0.95 & [1.57,1.99] & 0.74 \\
\end{tabular}
\caption{Best power-law regression fit (95\% confidence interval) for algorithm iteration runtime as a function of the number of non-zero values in the constraint matrix ($\fitexpr$).}
\label{tab:iteration_fit}
\end{table}
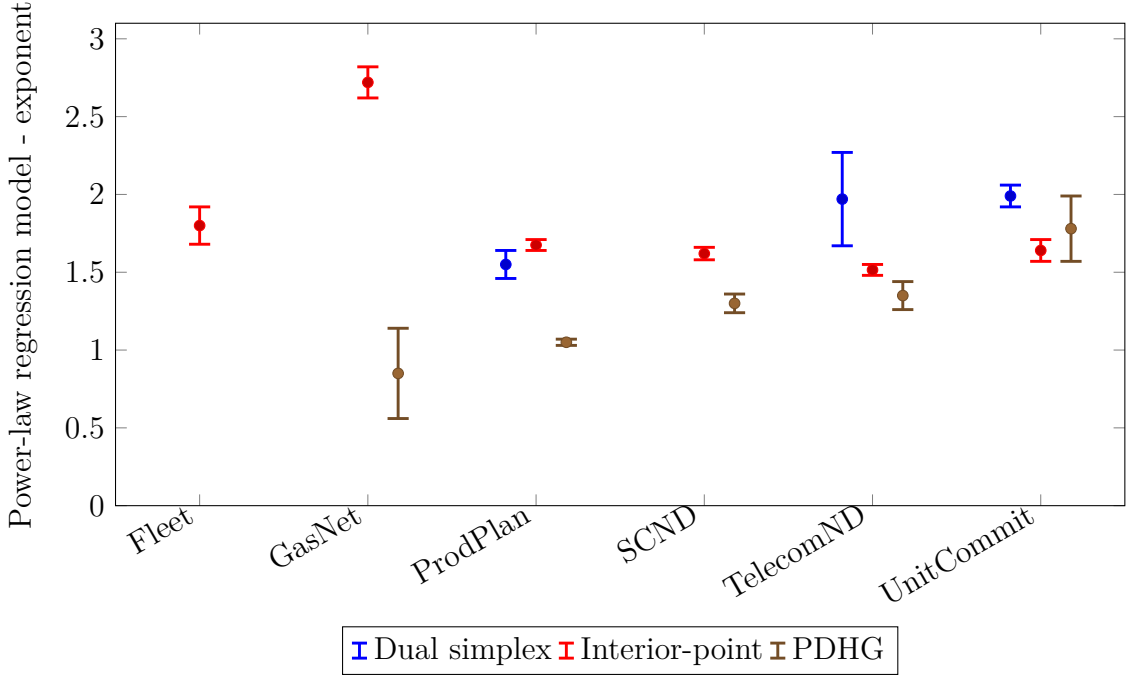
\begin{figure}[ht]
\centering
\pgfplotsset{
  ci plot/.style={
    only marks,
    mark=none,
    error bars/y dir=both,
    error bars/y explicit,
    error bars/error bar style={line width=1.1pt},
    error bars/error mark options={rotate=90, mark size=4pt, line width=1.1pt},
  }
}
\begin{tikzpicture}
\begin{axis}[
    width=15cm,
    height=8cm,
    ymin=0,
    ymax=3.1,
    xmin=0.5,
    xmax=6.5,
    xtick={1,2,3,4,5,6},
    xticklabels={Fleet,GasNet,ProdPlan,SCND,TelecomND,UnitCommit},
    xticklabel style={rotate=30, anchor=east},
    ylabel={Power-law regression model - exponent},
    legend style={
        at={(0.5,-0.25)},
        anchor=north,
        legend columns=-1
    },
    legend image code/.code={
      \draw[#1,line width=1.1pt] (0.15cm,-0.10cm) -- (0.15cm,0.10cm);
      \draw[#1,line width=1.1pt] (0.05cm,-0.10cm) -- (0.25cm,-0.10cm);
      \draw[#1,line width=1.1pt] (0.05cm, 0.10cm) -- (0.25cm, 0.10cm);
    },
]


\addplot+[ci plot]
coordinates {
    (2.82,1.55) +- (0,0.09)    
    (4.82,1.97) +- (0,0.30)    
    (5.82,1.99) +- (0,0.07)    
};
\addlegendentry{Dual simplex}

\addplot+[ci plot]
coordinates {
    (1.00,1.80)  +- (0,0.12)    
    (2.00,2.72)  +- (0,0.10)    
    (3.00,1.675) +- (0,0.035)   
    (4.00,1.62)  +- (0,0.04)    
    (5.00,1.515) +- (0,0.035)   
    (6.00,1.64)  +- (0,0.07)    
};
\addlegendentry{Interior-point}

\addplot+[ci plot]
coordinates {
    (2.18,0.85) +- (0,0.29)    
    (3.18,1.05) +- (0,0.02)    
    (4.18,1.30) +- (0,0.06)    
    (5.18,1.35) +- (0,0.09)    
    (6.18,1.78) +- (0,0.21)    
};
\addlegendentry{PDHG}

\end{axis}
\end{tikzpicture}
\caption{Best power-law regression fit (95\% confidence interval) for algorithm iteration runtime as a function of the number of non-zero values in the constraint matrix ($\fitexpr$).}
\label{fig:iteration_fit}
\end{figure}
The table gives three numbers for each algorithm/model pair: the 95\%
confidence interval for the predicted exponent in the power-law
regression fit ($b$ in ``$\fitexpr$''), and the $R^2$ result for that
fit (in log space).  Note that we do not include primal simplex here;
as noted earlier, dual simplex is much faster in general, so we allow
it to represent both.

The regression is performed over all of the instances that complete
their iterations, only reporting results when the algorithm solved at
least 20 of the 100 models in the set.  Dual simplex did not reach
this target for the Fleet, GasNet, and SCND sets, and PDHG did not
reach it for the Fleet set.

The regression model fits most cases quite well.  The interior-point
predictions have $R^2 > 0.92$ in all cases, which is generally
considered to be an excellent fit.  For dual simplex, $R^2 \geq 0.81$
in all cases where we have data.  For PDHG, $R^2 \geq 0.74$ for all
but the GasNet test set.

Among the algorithms considered, PDHG produces the smallest predicted
exponent in all but one case (UnitCommit).  Dual simplex produces the
largest exponent in all but one case.

\begin{figure}[htbp]
\centering
\begin{minipage}{0.48\textwidth}
\begin{tikzpicture}
\begin{loglogaxis}[
    xlabel={NZ},
    grid=none,
]
\addplot[
    only marks,
    mark=*,
] table[row sep=newline, col sep=comma] {
7.09858e+06,1782.6
7.25238e+06,1247.29
7.45301e+06,1968.87
7.58858e+06,3412.27
7.67697e+06,3025.34
7.78151e+06,2001.54
8.03185e+06,3203.31
8.13962e+06,4543.19
8.28517e+06,3443.35
8.48667e+06,3115.16
8.56315e+06,3694.23
8.79625e+06,7408.47
8.87881e+06,1782.42
9.04548e+06,6427.6
9.19985e+06,2518.41
9.325e+06,1902.22
9.5249e+06,2432.4
9.69357e+06,2938.78
9.8363e+06,2385.68
9.94968e+06,3305.42
1.01446e+07,5089.68
1.02991e+07,15432.9
1.04451e+07,2197.22
1.05685e+07,3628.32
1.07202e+07,4075.14
1.08201e+07,2798.03
1.1014e+07,2385.95
1.11538e+07,2827.23
1.12683e+07,2377.84
1.14348e+07,7015.46
1.16177e+07,3243.08
1.17665e+07,2503.22
1.18996e+07,4002.21
1.20591e+07,2140.81
1.22115e+07,12497.1
1.24045e+07,4179.11
1.25455e+07,4713.62
1.26534e+07,8184.39
1.28259e+07,4180.86
1.29143e+07,2482.16
1.31111e+07,3859.05
1.32694e+07,4619.55
1.33869e+07,3336.79
1.3589e+07,4681.96
1.36899e+07,3696.21
1.38404e+07,6212.99
1.40022e+07,4907.86
1.41788e+07,5153.71
1.43307e+07,9421.24
1.44975e+07,12263
1.4643e+07,4224.06
1.4804e+07,1796.25
1.48906e+07,10414.7
1.50429e+07,3159.39
1.52429e+07,3618.82
1.53468e+07,3861.62
1.55217e+07,5055.39
1.56415e+07,4113.29
1.57967e+07,4144.56
1.59616e+07,5554.58
1.60729e+07,9080.2
1.63394e+07,3775.25
1.64788e+07,4379.4
1.65884e+07,5556.65
1.67381e+07,5711.98
1.68354e+07,4661.94
1.70104e+07,4058.75
1.71986e+07,3511.16
1.72989e+07,4393.54
1.74455e+07,7109.03
1.76109e+07,4297.56
1.78403e+07,2742.28
1.7946e+07,3803.26
1.80146e+07,7949.19
1.82437e+07,8446.64
1.83616e+07,2772.14
1.85302e+07,3063.15
1.86936e+07,8429.35
1.87626e+07,8171.99
1.89822e+07,7180.26
1.90703e+07,6173.55
1.92644e+07,5667.26
1.93995e+07,5737.95
1.95867e+07,5693.72
1.97655e+07,5000.3
1.99214e+07,18811.4
1.9978e+07,5721
2.02214e+07,6754.02
2.0322e+07,14439.1
2.04664e+07,4172.68
2.06734e+07,3503.85
2.08468e+07,4170.97
2.09684e+07,8327.79
2.11206e+07,3924.57
2.12358e+07,9489.65
2.13733e+07,5412.39
2.15015e+07,3584.3
2.16731e+07,16103.5
2.18494e+07,3590.49
2.19213e+07,13894.5
};
\end{loglogaxis}
\end{tikzpicture}
\end{minipage}
\caption{Runtime versus non-zero count for PDHG iterations for GasNet.}
\label{fig:PDHGFleetGasNet}
\end{figure}
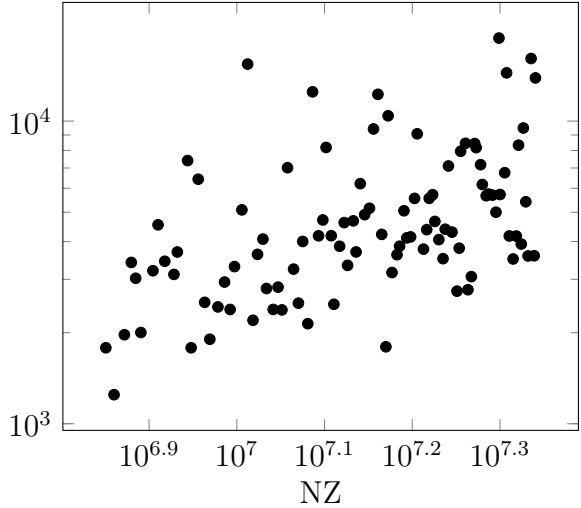

Let us now dig a bit deeper into the poor PDHG regression fit for the
GasNet test set.  One obvious question is whether this data
can be accurately modeled using a different function.
We tried fitting
polynomial, exponential, and log functions to the row, column,
and non-zero counts, and the quality of the model did not improve at all.
It is useful in situations like
this to also look at the raw data, which we show in
Figure~\ref{fig:PDHGFleetGasNet}.
PDHG runtimes are quite inconsistent, with no clear missed
opportunity to choose a better function to fit.

\subsection{Runtime Predictions by Algorithm Component}

As mentioned earlier, the three algorithms we consider for solving LPs
typically consist of several constituent parts.  The previous section
focused on the time required by the iterations performed by each
algorithm, but it is possible that iterations are not the
primary determinants of overall runtime.  We now look at predicted
growth rates for the full set of major algorithmic components.

We only consider the interior-point and PDHG algorithms here; the only
additional algorithm component for simplex beyond the simplex
iterations reported earlier is presolve, which will be reported here
for the other two algorithms.

\subsubsection{Interior-Point Components}

Recall that the interior-point method has four major components:
presolve, fill-reducing ordering, iterations, and (optionally)
crossover.  Table~\ref{tab:components_interiorpt} and
Figure~\ref{fig:components_interiorpt} show regression results for
these components for our six test sets.

\begin{table}[htbp]
\centering
\begin{tabular}{l|cc|cc|cc|cc}
  \textbf{} &
   \multicolumn{2}{c|}{Presolve} &
   \multicolumn{2}{c|}{Ordering} &
   \multicolumn{2}{c|}{Iterations} &
   \multicolumn{2}{c}{Crossover} \\
    &
   Exp (CI) & $R^2$ &
   Exp (CI) & $R^2$ &
   Exp (CI) & $R^2$ &
   Exp (CI) & $R^2$ \\ \hline
Fleet      & [1.00,1.03] & 1.00 & [2.20,2.39] & 0.96 & [1.68,1.92] & 0.92 & [1.75,1.97] & 0.94 \\
GasNet     & [1.12,1.14] & 1.00 & [2.24,2.47] & 0.95 & [2.62,2.82] & 0.97 & - & - \\
ProdPlan   & [1.09,1.10] & 1.00 & [1.48,1.54] & 0.99 & [1.64,1.71] & 0.99 & [0.99,1.03] & 0.99 \\
SCND       & [1.06,1.08] & 1.00 & [0.88,0.93] & 0.87 & [1.58,1.66] & 0.99 & [1.98,1.98] & 0.99 \\
TelecomND  & [1.07,1.10] & 1.00 & [0.04,0.31] & 0.06 & [1.48,1.55] & 0.99 & [1.57,1.89] & 0.82 \\
UnitCommit & [1.11,1.12] & 1.00 & [0.99,1.00] & 1.00 & [1.57,1.71] & 0.95 & [1.67,2.06] & 0.79 \\
\end{tabular}
\caption{Interior-point components: best power-law regression fit (95\% confidence interval) for component runtime against number of non-zero values in the constraint matrix ($\fitexpr$).}
\label{tab:components_interiorpt}
\end{table}

\begin{figure}[ht]
\centering

\pgfplotsset{
  ci plot/.style={
    only marks,
    mark=none,
    error bars/y dir=both,
    error bars/y explicit,
    error bars/error bar style={line width=1.1pt},
    error bars/error mark options={rotate=90, mark size=4pt, line width=1.1pt},
  }
}

\begin{tikzpicture}
\begin{axis}[
    width=15cm,
    height=8cm,
    ymin=0,
    ymax=3.1,
    enlarge x limits=0.15,
    symbolic x coords={Fleet,GasNet,ProdPlan,SCND,TelecomND,UnitCommit},
    xtick={Fleet,GasNet,ProdPlan,SCND,TelecomND,UnitCommit},
    xticklabel style={rotate=30, anchor=east},
    ylabel={Power-law regression model - exponent},
    legend style={
        at={(0.5,-0.25)},
        anchor=north,
        legend columns=-1
    },
    legend image code/.code={
      \draw[#1,line width=1.1pt] (0.15cm,-0.10cm) -- (0.15cm,0.10cm);
      \draw[#1,line width=1.1pt] (0.05cm,-0.10cm) -- (0.25cm,-0.10cm);
      \draw[#1,line width=1.1pt] (0.05cm, 0.10cm) -- (0.25cm, 0.10cm);
    },
]

\addplot+[ci plot,xshift=-9pt]
coordinates {
    (Fleet,1.015)      +- (0,0.015)
    (GasNet,1.13)      +- (0,0.01)
    (ProdPlan,1.095)   +- (0,0.005)
    (SCND,1.07)        +- (0,0.01)
    (TelecomND,1.085)  +- (0,0.015)
    (UnitCommit,1.115) +- (0,0.005)
};
\addlegendentry{Presolve}

\addplot+[ci plot,xshift=-3pt]
coordinates {
    (Fleet,2.295)      +- (0,0.095)
    (GasNet,2.355)     +- (0,0.115)
    (ProdPlan,1.51)    +- (0,0.03)
    (SCND,0.905)       +- (0,0.025)
    (TelecomND,0.175)  +- (0,0.135)
    (UnitCommit,0.995) +- (0,0.005)
};
\addlegendentry{Ordering}

\addplot+[ci plot,xshift=3pt]
coordinates {
    (Fleet,1.80)       +- (0,0.12)
    (GasNet,2.72)      +- (0,0.10)
    (ProdPlan,1.675)   +- (0,0.035)
    (SCND,1.62)        +- (0,0.04)
    (TelecomND,1.515)  +- (0,0.035)
    (UnitCommit,1.64)  +- (0,0.07)
};
\addlegendentry{Iterations}

\addplot+[ci plot,xshift=9pt]
coordinates {
    (Fleet,1.86)        +- (0,0.11)
    (ProdPlan,1.01)     +- (0,0.02)
    (SCND,1.98)         +- (0,0)
    (TelecomND,1.73)    +- (0,0.16)
    (UnitCommit,1.865)  +- (0,0.195)
};
\addlegendentry{Crossover}

\end{axis}
\end{tikzpicture}
\caption{Interior-point components: best power-law regression fit (95\% confidence interval) for component runtime against number of non-zero values in the constraint matrix ($\fitexpr$).}
\label{fig:components_interiorpt}
\end{figure}

The power-law regression models are generally quite good.  The model
achieves $R^2 = 1.00$ in all cases for presolve, $R^2 \geq 0.87$ in
all but one case for ordering, and $R^2 \geq 0.79$ in all cases for
crossover.  The table does not include data for crossover on the
GasNet set because crossover was unable to find a single basic
solution for any model in this set.

Presolve runtimes are nearly linear in the number of
non-zeros in $A$, which suggests that it will consume a diminishing
fraction of overall runtime as model sizes grow.

Note also that the exponents in the growth rate prediction are often
larger for crossover than they are for the interior-point iterations.
That represents a significant potential issue.  We will return to this
topic shortly.

For the Fleet set, ordering costs grow at a faster rate than interior-point
iteration costs.  This is surprising, but it is not outside the bounds
of what is known about ordering algorithms.  Perhaps somewhat
ironically, ordering heuristics are meant to reduce the cost of the
subsequent sparse factorization, but there are no guarantees that
their cost is smaller than that of the factorization.

\subsubsection{PDHG Components}

\begin{figure}[ht]
\centering
\pgfplotsset{
  ci plot/.style={
    only marks,
    mark=none,
    error bars/y dir=both,
    error bars/y explicit,
    error bars/error bar style={line width=1.1pt},
    error bars/error mark options={rotate=90, mark size=4pt, line width=1.1pt},
  }
}

\begin{tikzpicture}
\begin{axis}[
    width=15cm,
    height=8cm,
    ymin=0,
    ymax=2.8,
    enlarge x limits=0.15,
    symbolic x coords={Fleet,GasNet,ProdPlan,SCND,TelecomND,UnitCommit},
    xtick={Fleet,GasNet,ProdPlan,SCND,TelecomND,UnitCommit},
    xticklabel style={rotate=30, anchor=east},
    ylabel={Power-law regression model - exponent},
    legend style={
        at={(0.5,-0.25)},
        anchor=north,
        legend columns=-1
    },
    legend image code/.code={
      \draw[#1,line width=1.1pt] (0.15cm,-0.10cm) -- (0.15cm,0.10cm);
      \draw[#1,line width=1.1pt] (0.05cm,-0.10cm) -- (0.25cm,-0.10cm);
      \draw[#1,line width=1.1pt] (0.05cm, 0.10cm) -- (0.25cm, 0.10cm);
    },
]

\addplot+[ci plot,xshift=-6pt]
coordinates {
    (Fleet,1.02)       +- (0,0.01)
    (GasNet,1.135)     +- (0,0.005)
    (ProdPlan,1.10)    +- (0,0.00)
    (SCND,1.07)        +- (0,0.01)
    (TelecomND,1.095)  +- (0,0.015)
    (UnitCommit,1.11)  +- (0,0.01)
};
\addlegendentry{Presolve}

\addplot+[ci plot]
coordinates {
    (GasNet,0.85)       +- (0,0.29)
    (ProdPlan,1.05)     +- (0,0.02)
    (SCND,1.30)         +- (0,0.06)
    (TelecomND,1.35)    +- (0,0.09)
    (UnitCommit,1.78)   +- (0,0.21)
};
\addlegendentry{Iterations}

\addplot+[ci plot,xshift=6pt,draw=black,mark=*, mark options={fill=black}]
coordinates {
    (ProdPlan,1.145)     +- (0,0.015)
    (SCND,2.295)         +- (0,0.235)
    (TelecomND,1.835)    +- (0,0.125)
    (UnitCommit,2.03)    +- (0,0.23)
};
\addlegendentry{Crossover}

\end{axis}
\end{tikzpicture}
\caption{PDHG components: best power-law regression fit (95\% confidence interval) for component runtime against number of non-zero values in the constraint matrix ($\fitexpr$).}
\label{fig:components_pdhg}
\end{figure}
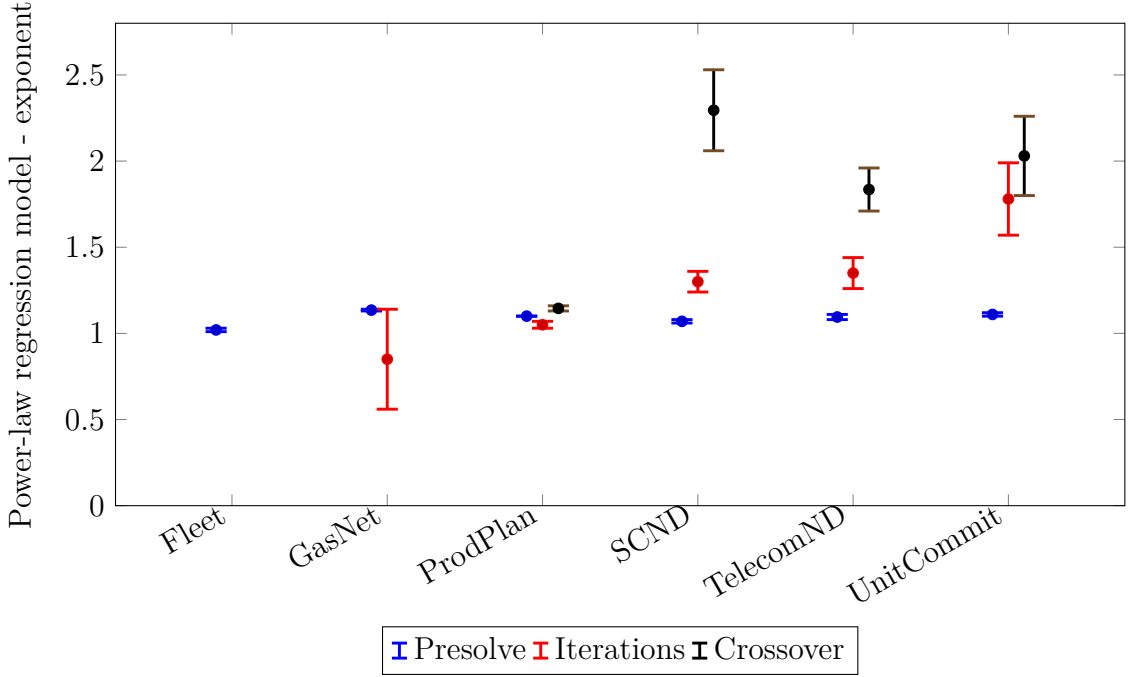

Recall that PDHG has three major components: presolve, PDHG iterations, and (optionally) crossover.
Table~\ref{tab:components_pdhg} and Figure~\ref{fig:components_pdhg} show
regression results for these components.
\begin{table}[htbp]
\centering
\begin{tabular}{l|cc|cc|cc}
  \textbf{} &
   \multicolumn{2}{c|}{Presolve} &
   \multicolumn{2}{c|}{Iterations} &
   \multicolumn{2}{c}{Crossover} \\
    &
   Exp (CI) & $R^2$ &
   Exp (CI) & $R^2$ &
   Exp (CI) & $R^2$ \\ \hline
Fleet      & [1.01,1.03] & 1.00 & - & - & - & - \\
GasNet     & [1.13,1.14] & 1.00 & [0.56,1.14] & 0.26 & - & - \\
ProdPlan   & [1.10,1.10] & 1.00 & [1.03,1.07] & 0.99 & [1.13,1.16] & 1.00 \\
SCND       & [1.06,1.08] & 1.00 & [1.24,1.36] & 0.95 & [2.06,2.53] & 0.86 \\
TelecomND  & [1.08,1.11] & 0.99 & [1.26,1.44] & 0.91 & [1.71,1.96] & 0.90 \\
UnitCommit & [1.10,1.12] & 1.00 & [1.57,1.99] & 0.74 & [1.80,2.26] & 0.76 \\
\end{tabular}
\caption{PDHG components: best power-law regression fit (95\% confidence interval) for component runtime against number of non-zero values in the constraint matrix ($\fitexpr$).}
\label{tab:components_pdhg}
\end{table}

A few things are noteworthy in these results.  The first is that some
results are missing.  PDHG was only able to solve 16 of the 100 models
in our Fleet set, so it didn't meet our threshold of 20.  Similarly,
crossover results are missing for Fleet and Gasnet; we were only able
to find basic solutions for 8 of the 100 models in the Fleet set
within our 50,000-second time limit, and we did not find a single
basic solution for GasNet.

Next, note that predicted growth rates for crossover runtimes are
larger than those for PDHG in every case.  When combined with the fact
that crossover is inherently sequential while PDHG can make excellent
use of parallel computing, this suggests that without some new ideas,
crossover is likely to dominate the overall cost of finding a basic
solution with PDHG.  Crossover is of course optional, but it is
currently the best method available to remedy the low accuracy of PDHG
solutions.

Finally, when comparing this data against the data for the
interior-point solver in Table~\ref{tab:components_interiorpt}, the
crossover growth rates appear to be consistently larger for PDHG.
That partially fits our expectations, since the final PDHG iterates
are much less accurate (although we were expecting this to be a
constant-factor effect, rather than an asymptotic effect).

\subsection{Absolute Performance Comparison}

While understanding asymptotic behavior is important, lower-order
terms definitely matter as well.  Our final test compares performance
of the interior-point and PDHG solvers on each of our test sets using
measured runtimes on our test system with a modern CPU and GPU (an EPYC 9575F
and an Nvidia RTX Pro 6000, respectively).  For each test set, we show
results for the three largest models that solve to optimality within a
10,000-second time limit.  We include three results because
performance varied enough from instance to instance that a single data
point could give the wrong impression.

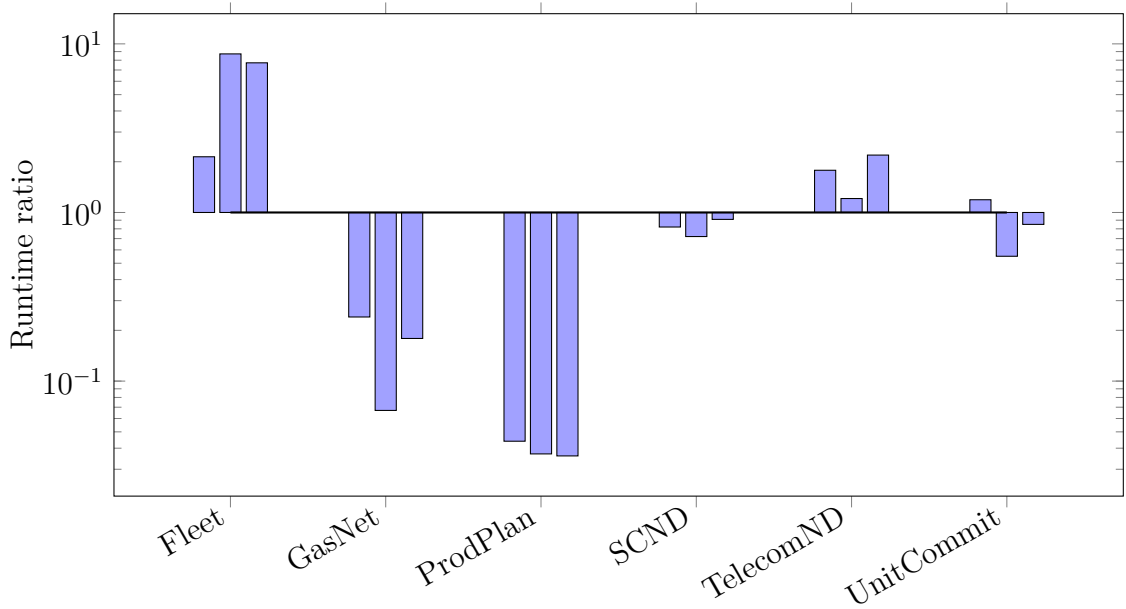
\begin{figure}[bthp]
\centering
\begin{tikzpicture}
\begin{axis}[
    ybar,
    ymode=log,
    bar width=8pt,
    width=15cm,
    height=8cm,
    enlarge x limits=0.15,
    symbolic x coords={Fleet,GasNet,ProdPlan,SCND,TelecomND,UnitCommit},
    xtick=data,
    xticklabel style={rotate=30, anchor=east},
    ylabel={Runtime ratio},
    legend style={
        at={(0.5,-0.18)},
        anchor=north,
        legend columns=-1
    },
]

\addplot[fill=blue!37] coordinates {
    (Fleet,2.14) (GasNet,0.24) (ProdPlan,0.044)
    (SCND,0.82) (TelecomND,1.78) (UnitCommit,1.19)
};
\addplot[fill=blue!37] coordinates {
    (Fleet,8.72) (GasNet,0.067) (ProdPlan,0.037)
    (SCND,0.72) (TelecomND,1.21) (UnitCommit,0.55)
};
\addplot[fill=blue!37] coordinates {
    (Fleet,7.72) (GasNet,0.179) (ProdPlan,0.036)
    (SCND,0.91) (TelecomND,2.19) (UnitCommit,0.85)
};

\draw[black, thick] (axis cs:Fleet,1) -- (axis cs:UnitCommit,1);

\end{axis}
\end{tikzpicture}
\caption{Ratio of PDHG runtime to interior-point runtime (crossover disabled) for the three largest
instances that both methods can solve (interior-point method on the CPU, PDHG on the GPU).}
\label{fig:absolute_nocrossover}
\end{figure}

Figure~\ref{fig:absolute_nocrossover} shows the ratio of PDHG runtime
to interior-point runtime without crossover.
The interior-point solver provides better performance on one test set
(Fleet), and PDHG provides better performance on two (GasNet and
ProdPlan).  They are basically tied on the other three.  Recall that
asymptotic growth rates favor PDHG on most of these sets, suggesting
that other factors still play a significant role at these model sizes.
Recall that this is an asymmetric comparison, since PDHG will
typically give errors that are 4-5 orders of magnitude larger with
these termination criteria.

\begin{figure}[bthp]
\centering
\begin{tikzpicture}
\begin{axis}[
    ybar,
    ymode=log,
    bar width=8pt,
    width=15cm,
    height=8cm,
    enlarge x limits=0.15,
    symbolic x coords={Fleet,GasNet,ProdPlan,SCND,TelecomND,UnitCommit},
    xtick=data,
    xticklabel style={rotate=30, anchor=east},
    ylabel={Runtime ratio},
    legend style={
        at={(0.5,-0.18)},
        anchor=north,
        legend columns=-1
    },
]

\addplot[fill=blue!37] coordinates {
    (Fleet,30) (GasNet,nan) (ProdPlan,0.0702)
    (SCND,7.44) (TelecomND,1.82) (UnitCommit,0.816)
};
\addplot[fill=blue!37] coordinates {
    (Fleet,30) (GasNet,nan) (ProdPlan,0.0653)
    (SCND,12.44) (TelecomND,1.65) (UnitCommit,0.726)
};
\addplot[fill=blue!37] coordinates {
    (Fleet,30) (GasNet,nan) (ProdPlan,0.0572)
    (SCND,7.68) (TelecomND,1.09) (UnitCommit,0.745)
};

\draw[black, thick] (axis cs:Fleet,1) -- (axis cs:UnitCommit,1);

\end{axis}
\end{tikzpicture}
\caption{Ratio of PDHG runtime to interior-point runtime (crossover enabled) for the three largest
  instances that both methods can solve (interior-point method and crossover on the CPU, PDHG on the GPU).}
\label{fig:absolute_crossover}
\end{figure}
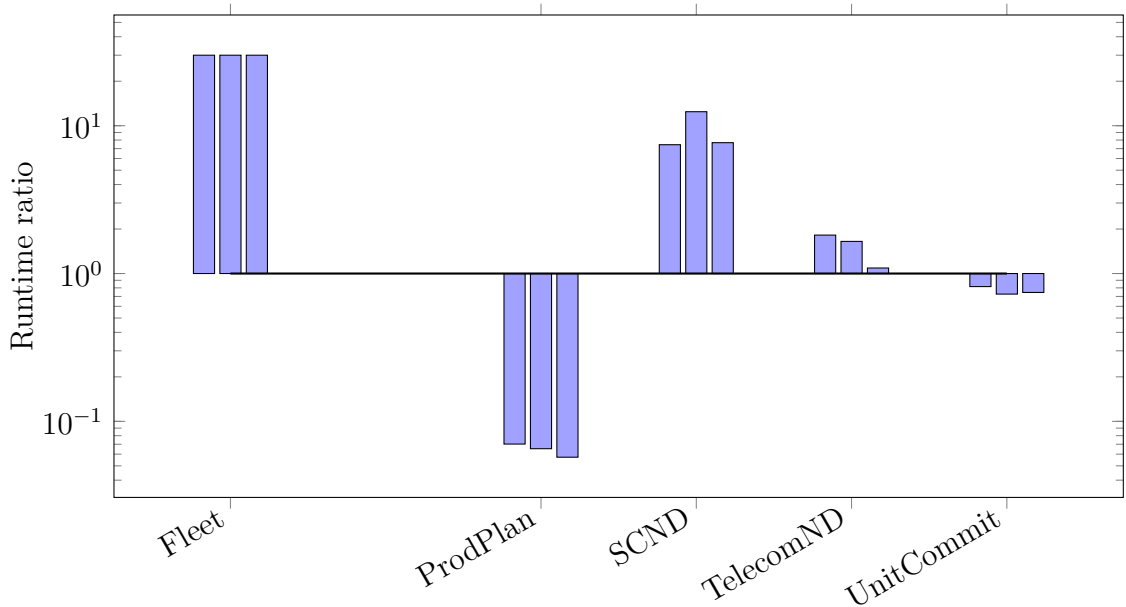

Figure~\ref{fig:absolute_crossover} shows the same comparison, but with crossover enabled.
Both algorithms now produce basic solutions, so differences in
solution quality are no longer an issue.  We capped the ratio at 30 for
the Fleet test set, because crossover rarely converged for these models
when using PDHG.  Data is missing for the GasNet models because crossover was
unable to produce basic solutions for any instances in our test set.
Among the five test sets for which both
methods produced a basic solution, the interior-point solver is
significantly faster on two (Fleet and SCND), significantly
slower on one (ProdPlan), and results are similar for two (TelecomND
and UnitCommit).

\begin{figure}[bthp]
\centering
\begin{tikzpicture}
\begin{axis}[
    ybar,
    bar width=8pt,
    width=15cm,
    height=8cm,
    enlarge x limits=0.18,
    symbolic x coords={Fleet,ProdPlan,SCND,TelecomND,UnitCommit},
    xtick=data,
    xticklabel style={rotate=30, anchor=east},
    ylabel={Fraction of runtime in crossover (\%)},
    ymin=0,
    ymax=100,
    ytick={0,20,40,60,80,100},
    legend style={
        at={(0.5,-0.35)},
        anchor=south,
        legend columns=-1
    },
]

\addplot[fill=blue!37, draw=black, bar shift=-5pt, forget plot] coordinates {
    (Fleet,6.8)
    (ProdPlan,0.9)
    (SCND,46)
    (TelecomND,42)
    (UnitCommit,25)
};

\addplot[fill=orange!60, draw=black, bar shift=5pt, forget plot] coordinates {
    (Fleet,nan)
    (ProdPlan,34)
    (SCND,95)
    (TelecomND,45)
    (UnitCommit,18)
};

\addlegendimage{area legend, fill=blue!37, draw=black}
\addlegendentry{Interior-point}
\addlegendimage{area legend, fill=orange!60, draw=black}
\addlegendentry{PDHG}

\end{axis}
\end{tikzpicture}
\caption{Fraction of total runtime spent in crossover for interior-point and PDHG methods.}
\label{fig:crossover_fraction}
\end{figure}
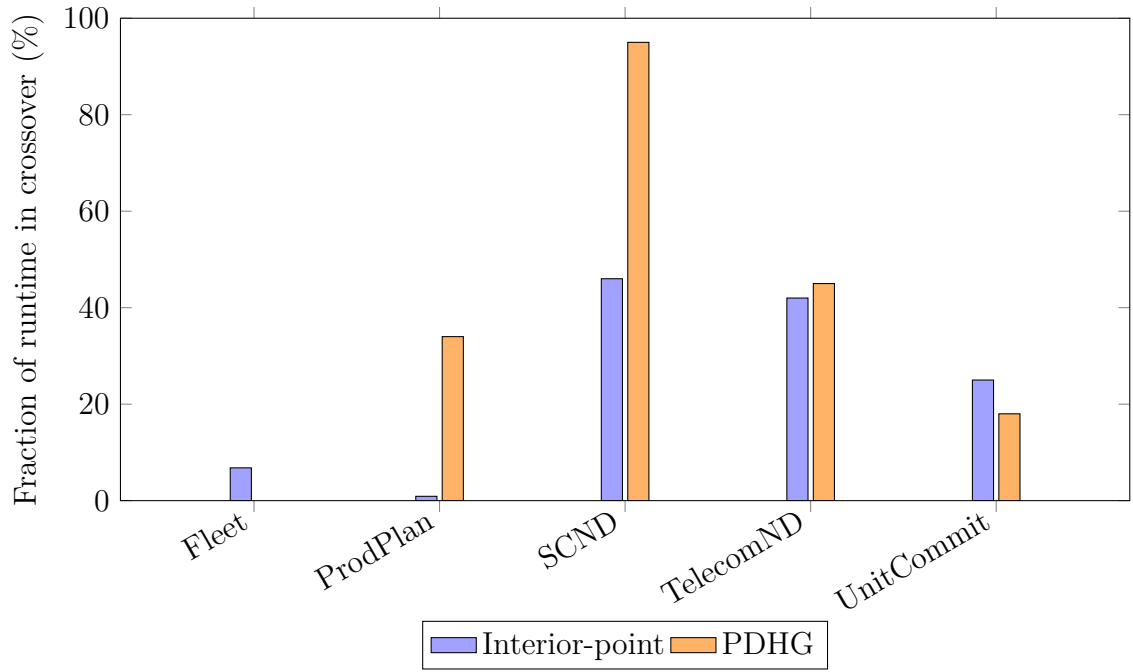

Figure~\ref{fig:crossover_fraction} shows the fraction of overall runtime
spent in crosover, averaged over the three models from each test set
considered in the previous figure.
Crossover times do not dominate total runtime for the interior-point
solver, but they do represent a significant fraction in several cases.
Crossover for PDHG consistently consumes a larger fraction of overall
time than it does for interior-point solver, but this is not
universally true.  The lower accuracy from PDHG clearly has an impact,
but not always.

\section{Discussion}

The results presented in this paper suggest that there is a
fundamental tradeoff between speed and accuracy for LP algorithms.
Generally speaking, PDHG has the lowest predicted asymptotic runtime,
followed by the interior-point method, followed by simplex/crossover.
In contrast, PDHG provides the lowest accuracy, followed by the
interior-point method, followed by the basic solutions computed by
simplex and crossover.

A number of ideas have been proposed for moving around in this space,
including {\em feasibility polishing\/}~\cite{pdlp25}, new
PDHG-inspired crossover algorithms~\cite{pdhgcrossover24}, concurrent
crossover~\cite{concurrentcrossover25}, PDHG corner
pushes~\cite{pdhgcorner25}, and hybrid PDHG/interior-point
methods~\cite{pdhghybrid26}.  While these may help in specific
instances, so far none seem to alter the fundamental calculus
noted above.

This brings us back to the question of how much accuracy is enough.
The answer will ultimately come from the consumers of these solutions.
There have been a few reports of success in using low-accuracy
solutions in novel ways (for example, to guide a MIP heuristic to a
high-quality, high-accuracy integer solution~\cite{pdhgheuristic25}).
At this point, it is not clear whether this will be the exception
or the rule.

Another question that arises from the results in this paper is whether
synthetic models and data can provide useful insights into the
behavior of practical optimization models.  As we noted, diving into
the models made us fairly comfortable with their realism.  The data is
the tougher part.  Making progress on this will require some objective
measure of model realism, which appears to be a difficult task.

\section{Conclusions}

This paper looked at how runtimes grow with model size for the three
most powerful algorithms for solving large linear-programming
problems: simplex, interior-point, and PDHG.  We considered six test
sets that were inspired by practical applications that require the
solution of very large LPs.  To be able to consider a wide range of
input sizes, we exploited the recent ability of large language models
to generate credible optimization models, along with synthetic data to
feed these models.

Our results show clear differences in runtime growth rates of the
different algorithms.  Some observed behaviors were different for the
different model types, while others were quite consistent.  In
general, the methods that produced the most accurate solutions gave
the worst predicted asymptotic behavior.

\section*{Acknowledgements}

We are grateful to Everett Dutton for insightful discussions that
greatly improved the presentation of the results.

\bibliographystyle{plainnat}
\bibliography{references}

\begin{thebibliography}{23}
\providecommand{\natexlab}[1]{#1}
\providecommand{\url}[1]{\texttt{#1}}
\expandafter\ifx\csname urlstyle\endcsname\relax
  \providecommand{\doi}[1]{doi: #1}\else
  \providecommand{\doi}{doi: \begingroup \urlstyle{rm}\Url}\fi

\bibitem[Andersen and Andersen(1996)]{presolve1995}
Erling~D. Andersen and Knud~D. Andersen.
\newblock Presolving in linear programming.
\newblock \emph{Mathematical Programming}, 71:\penalty0 221--245, 1996.

\bibitem[Applegate et~al.(2022)Applegate, D{\'i}az, Hinder, Lu, Lubin,
  O’Donoghue, and Schudy]{pdlp22}
David Applegate, Mateo D{\'i}az, Oliver Hinder, Haihao Lu, Miles Lubin, Brendan
  O’Donoghue, and Warren Schudy.
\newblock Practical large-scale linear programming using primal-dual hybrid
  gradient.
\newblock \emph{Advances in Neural Information Processing Systems}, 34\penalty0
  (2021):\penalty0 20243–--20257, 2022.

\bibitem[Applegate et~al.(2025)Applegate, D{\'i}az, Hinder, Lu, Lubin,
  O’Donoghue, and Schudy]{pdlp25}
David Applegate, Mateo D{\'i}az, Oliver Hinder, Haihao Lu, Miles Lubin, Brendan
  O’Donoghue, and Warren Schudy.
\newblock {PDLP}: a practical first-order method for large-scale linear
  programming.
\newblock \emph{arXiv preprint arXiv:2501.07018}, 2025.
\newblock URL \url{https://arxiv.org/abs/2501.07018}.

\bibitem[Chambolle and Pock(2011)]{pdhg11}
Antonin Chambolle and Thomas Pock.
\newblock A first-order primal-dual algorithm for convex problems with
  applications to imaging.
\newblock \emph{Journal of Mathematical Imaging and Vision}, 40:\penalty0
  120--145, 2011.

\bibitem[Chen et~al.(2024)Chen, Sun, Yuan, Zhang, and Zhao]{hprlp24}
Kaihuang Chen, Defeng Sun, Yancheng Yuan, Guojun Zhang, and Xinyuan Zhao.
\newblock {HPR-LP}: An implementation of an {HPR} method for solving linear
  programming.
\newblock \emph{arXiv preprint arXiv:2408.12179}, 2024.
\newblock URL \url{https://arxiv.org/abs/2408.12179}.

\bibitem[Dantzig(1951)]{simplex51}
George~B. Dantzig.
\newblock Applications of the simplex method to a transportation problem.
\newblock In \emph{Activity Analysis of Production and Allocation}, pages
  359--373, 1951.

\bibitem[Demmel(1997)]{laplacian}
J.W. Demmel.
\newblock \emph{Applied Numerical Linear Algebra}.
\newblock SIAM, 1997.

\bibitem[Draper and Smith(1998)]{regression}
N.R. Draper and H.~Smith.
\newblock \emph{Applied Regression Analysis}.
\newblock Wiley, 1998.

\bibitem[George and Liu(1981)]{georgeliu81}
A.~George and J.W.H. Liu.
\newblock \emph{Computer Solution of Large Sparse Positive Definite Systems}.
\newblock Prentice Hall, 1981.

\bibitem[Kempke and Koch(2025)]{pdhgheuristic25}
Nils-Christian Kempke and Thorsten Koch.
\newblock Fix-and-propage heuristics using low-precision first-order {LP}
  solutions for large-scale mixed-integer linear optimization.
\newblock \emph{arXiv preprint arXiv:2503.10344}, 2025.
\newblock URL \url{https://arxiv.org/abs/2503.10344}.

\bibitem[Liu and Lu(2024)]{pdhgcrossover24}
Tianhao Liu and Haihao Lu.
\newblock A new crossover algorithm for {LP} inspired by the spiral dynamic of
  {PDHG}.
\newblock \emph{arXiv preprint arXiv:2409.14715}, 2024.
\newblock URL \url{https://arxiv.org/abs/2409.14715}.

\bibitem[Lu and Yang(2023)]{cupdlp23}
Haihao Lu and Jinwen Yang.
\newblock {cuPDLP.jl}: A {GPU} implementation of restarted primal-dual hybrid
  gradient for linear programming in {J}ulia.
\newblock \emph{arXiv preprint arXiv:2311.12180}, 2023.
\newblock URL \url{https://arxiv.org/abs/2311.12180}.

\bibitem[Lu et~al.(2025)Lu, Peng, and Yang]{cupdlp+25}
Haihao Lu, Zedong Peng, and Jinwen Yang.
\newblock {cuPDLPx}: A further enhanced {GPU}-based first-order solver for
  linear programming.
\newblock \emph{arXiv preprint arXiv:2507.14051}, 2025.
\newblock URL \url{https://arxiv.org/abs/2507.14051}.

\bibitem[Megiddo(1991)]{crossover91}
Nimrod Megiddo.
\newblock On finding primal- and dual-optimal bases.
\newblock \emph{ORSA Journal on Computing}, 3\penalty0 (1):\penalty0 63--65,
  1991.

\bibitem[Mittelmann(2025)]{hans25}
Hans Mittelmann.
\newblock Decision tree for optimization software.
\newblock \url{http://plato.asu.edu/guide.html}, 2025.

\bibitem[{OpenAI}(2026)]{chatgpt2026}
{OpenAI}.
\newblock Chat{GPT}, 2026.
\newblock URL \url{https://chatgpt.com}.

\bibitem[Optimization(2025)]{gurobot}
Gurobi Optimization.
\newblock Gurobot, 2025.
\newblock URL \url{https://www.gurobi.com/solutions/gurobot}.

\bibitem[Rothberg(2025{\natexlab{a}})]{concurrentcrossover25}
Edward Rothberg.
\newblock Concurrent crossover for {PDHG}.
\newblock \emph{arXiv preprint arXiv:2510.24429}, 2025{\natexlab{a}}.
\newblock URL \url{https://arxiv.org/abs/2510.24429}.

\bibitem[Rothberg(2025{\natexlab{b}})]{pdhgcorner25}
Edward Rothberg.
\newblock Backing {PDHG} into a corner.
\newblock \emph{arXiv preprint arXiv:2511.13894}, 2025{\natexlab{b}}.
\newblock URL \url{https://arxiv.org/abs/2511.13894}.

\bibitem[Rothberg(2026)]{pdhghybrid26}
Edward Rothberg.
\newblock Hybridizing {PDHG} and interior-point methods.
\newblock \emph{arXiv preprint arXiv:2603.03150}, 2026.
\newblock URL \url{https://arxiv.org/abs/2603.03150}.

\bibitem[Stoian et~al.(2022)Stoian, Guinchiglia, and Lukasiewicz]{synthetic25}
Mihaela~Catalina Stoian, Eleonora Guinchiglia, and Thomas Lukasiewicz.
\newblock A survey on tabular data generation: Utility, alignment, fidelity,
  privacy, and beyond.
\newblock \emph{arXiv preprint arXiv:2503.05954v1}, 2022.
\newblock URL \url{https://arxiv.org/html/2503.05954v1}.

\bibitem[Wright(1997)]{wright97}
Stephen Wright.
\newblock \emph{Primal-dual interior-point methods}.
\newblock SIAM, 1997.

\bibitem[Zhang et~al.(2025)Zhang, Chen, Zope, Barbalho, Mellou, Molinaro,
  Kulkarni, Menache, and Li]{zhang2025optimind}
Xinzhi Zhang, Zeyi Chen, Humishka Zope, Hugo Barbalho, Konstantina Mellou,
  Marco Molinaro, Janardhan Kulkarni, Ishai Menache, and Sirui Li.
\newblock {OptiMind}: Teaching {LLM}s to think like optimization experts.
\newblock \emph{arXiv preprint arXiv:2509.22979}, 2025.

\end{thebibliography}

\end{document}